\documentclass[11pt]{article}

\usepackage{geometry}                
\usepackage{graphicx}
\usepackage{amssymb}
\usepackage{epstopdf}
\usepackage{amsmath,amsthm,amscd,amssymb}
\usepackage{latexsym}
\numberwithin{equation}{section}

\theoremstyle{plain}
\newtheorem{theorem}{Theorem}[section]
\newtheorem{lemma}[theorem]{Lemma}

\theoremstyle{definition}

\theoremstyle{remark}
\newtheorem{remark}[theorem]{Remark}

\newtheorem{case[theorem]}{Case}

\newcommand{\ds}{\displaystyle}
\newcommand{\supp}{\mathrm{supp}\,}

\title{A universal Stein-Tomas restriction estimate for measures in three dimensions}
\author{Alex Iosevich and Svetlana Roudenko}
\date{}

\begin{document}

\maketitle

\begin{abstract}
We study restriction estimates in ${\Bbb R}^3$ for surfaces given as
graphs of $W^1_1({\Bbb R}^2)$ (integrable gradient) functions. We
obtain a ``universal"
$$
L^2(\mu) \to L^4({\Bbb R}^3, L^2(SO(3)) \,)
$$
estimate for the extension operator $f \to \widehat{f \mu}$ in three
dimensions. We also prove that the three dimensional estimate holds
for any Frostman measure supported on a compact set of Hausdorff
dimension greater than two. The approach is geometric and is
influenced by a connection with the Falconer distance problem.
\end{abstract}

\tableofcontents

\section{Introduction}

The classical Stein-Tomas restriction theorem says that if $\mu$ is
the Lebesgue measure on $S^{d-1}$, the unit sphere in ${\Bbb R}^d$,
or, more generally, on a smooth convex surface with everywhere
non-vanishing curvature, then
\begin{equation}
 \label{steintomas}
{\Vert \widehat{f \mu} \Vert}_{L^{\frac{2(d+1)}{d-1}}({\Bbb R}^d)}
\lesssim {\Vert f \Vert}_{L^2(S^{d-1})}.
\end{equation}

It is shown in \cite{IL00} (see also \cite{I99}) that if the
Gaussian curvature is allowed to vanish, (\ref{steintomas}) does not
hold. Nevertheless, there is hope of obtaining (\ref{steintomas})
for all reasonable surfaces by modifying the surface carried measure
in some universal way. For example, if $\mu_0$ is the Lebesgue
measure on a convex compact smooth surface $\Gamma$, finite type in
the sense that the order of contact with every tangent line is
finite, and $d\mu(x)=K^{\frac{1}{d+1}}(x) d\mu_0(x)$, then one can
check using standard techniques that the estimate (\ref{steintomas})
holds. The situation becomes much more complicated if the Gaussian
curvature is allowed to vanish to infinite order. Carbery, Kenig and
Ziesler \cite{CKZ06} recently proved (\ref{steintomas}) for suitably
weighted measures on surfaces of revolution in three dimensions
under some quantitative assumption on the graphing function.

In \cite{BIT01}, the authors took a different point of view. Instead
of imposing a fixed measure on the family of surfaces, they
considered mixed norm restriction theorems corresponding to convex
curves under rotations. The approach was heavily tied to the average
decay estimates of the Fourier transform of the Lebesgue measure on
convex curves, due to Podkorytov, which made the convexity
assumption difficult to by-pass. In this paper, we take a geometric
point of view which allows us to consider a much more general
collection of surfaces. Our main result is the following.

\begin{theorem}
 \label{main}
Let $\mu$ be the Frostman measure on a compact $two$-dimensional
surface $S$ in ${\Bbb R}^3$ given as the graph of a $W^1_1({\Bbb
R}^2)$ function. Recall that $W^1_1({\Bbb R}^2)$ is the class of
functions in two variables whose gradient is in $L^1({\Bbb R}^2)$.

Given $\theta \in SO(3)$, $d \ge 3$, the special orthogonal group,
define the random measure $d\mu_{\theta}$ by the formula
$$
\int g(x) \,d\mu_{\theta}(x)=\int f(\theta x) \,d\mu(x).
$$
Then,
\begin{equation}
 \label{mainest}
{\left( \int {\left| \int_{SO(3)} {|\widehat{f \mu_{\theta}(x)}|}^2
\,dH(\theta) \right|}^2 dx \right)}^{\frac{1}{4}} \lesssim
{||f||}_{L^2(d\mu)},
\end{equation}
where $dH(\theta)$ is the Haar measure on $SO(3)$.

Moreover, the same estimate holds if $\mu$ is the Frostman measure
on any compact subset of ${\mathbb R}^3$ of Hausdorff dimension
greater than two.
\end{theorem}

\begin{remark}
The condition we need to impose on the measure $\mu$ in order for
the conclusion of Theorem \ref{main} to hold is that
\begin{equation}
 \label{E:mu-epsilon}
\mu \times \mu \{(x,y): 1 \leq |x-y| \leq 1+\epsilon \} \lesssim
\epsilon.
\end{equation}
This condition holds for Lipschitz surfaces, but it also holds for
many measures supported on sets that are far from rectifiable in any
sense. For example, consider a sequence of positive integers
$\{q_i\}$ such that $q_1=2$ and $q_{i+1}>q_i^i$. Let $E_q$ denote
the $q^{-\frac{d}{s}}$, $0<s<d$, neighborhood of $q^{-1}({\Bbb Z}^d
\cap {[0,q]}^d)$. Let $\ds E_s=\cap_{i=1}^{\infty} E_{q_i}$. By
standard geometric measure theory (see e.g. \cite{Fal85}), the
Hausdorff dimension of $E_s$ is $s$. Let $s=d-1$. One can check by a
direct calculation that (\ref{E:mu-epsilon}) holds.
\end{remark}

\vskip.125in

\subsection{Acknowledgements:}
The author wishes to thank Michael Loss of Georgia Institute of
Technology for a helpful suggestion regarding the regularity
assumptions in the main result. S.R. was partially supported by the
NSF grant DMS-0531337.

\vskip.125in

\section{Reduction to the key geometric estimate}


Let
$$
T(f,g)(x)=\int f\mu_{\theta}*g\mu_{\theta}(x) \,dH_d(\theta),
$$
where $H$ is the Haar (probability) measure on $SO(d)$.

On one hand,
$$
{||T(f,g)||}_{L^1({\Bbb R}^d)} \lesssim {||f||}_{L^1(\mu)}
\cdot {||g||}_{L^1(\mu)},
$$
since convolution of two $L^1$ functions is in $L^1$ by Fubini.

On the other hand,
$$
{||T(f,g)||}_{L^{\infty}({\Bbb R}^d)} \lesssim
{||f||}_{L^{\infty}(\mu)} \cdot {||g||}_{L^{\infty}(\mu)} \cdot
\sup_{x \in {\Bbb R}^d} \left| \int \mu_{\theta}*\mu_{\theta}(x)\,
dx \right|.
$$
It follows by interpolation and setting $f=g$ that if
\begin{equation}
 \label{convolutionbound}
\sup_{x \in {\Bbb R}^d} \left| \int \mu_{\theta}*\mu_{\theta}(x) \,
dH(\theta) \right| \lesssim 1,
\end{equation}
then
\begin{equation}
 \label{basicestimate}
{\left( \int {\left( \int {|\widehat{f \mu_{\theta}}(x)|}^2
dH_d(\theta) \right)}^2 dx \right)}^{\frac{1}{4}} \lesssim
{||f||}_{L^2(\mu)}.
\end{equation}
This reduces matters to the study of (\ref{convolutionbound}) and
this is what the remainder of the paper is about. Since
$$
\supp(\mu_{\theta}*\mu_{\theta}) \subset \supp(\mu_{\theta})+
\supp(\mu_{\theta}),
$$
we can easily arrange to take a supremum over $x$ away from a fixed
neighborhood of the origin. This is precisely what we shall do in
the sequel.

\vskip.125in

\section{Proof of Theorem \ref{main}}

By the Fourier inversion formula
$$ \int \mu_{\theta}*\mu_{\theta}(x) \,dH(\theta)$$
$$
=\int \int e^{2 \pi i \theta x \cdot \xi} {|\widehat{\mu}(\xi)|}^2
\, dH(\theta) \,d\xi
$$
$$
=\int \int_{S^2} \int_{\{\theta \in SO(3): \theta
\frac{x}{|x|}=\omega\}} e^{2 \pi i \omega \cdot |x| \xi} \,
dH(\theta) \, d\omega {|\widehat{\mu}(\xi)|}^2 \, d\xi
$$
\begin{equation}
 \label{stationary3}
=\int \int_{S^{2}} e^{2 \pi i \omega \cdot |x| \xi} h(\omega) \,
d\omega {|\widehat{\mu}(\xi)|}^2 \, d\xi,
\end{equation}
where
$$
h(\omega)=H \left(\left\{\theta \in SO(3): \theta
\frac{x}{|x|}=\omega \right\} \right).
$$
\begin{lemma}
 \label{smooth3}
The function $h(\omega)$ is constant.
\end{lemma}

The proof is immediate since the Haar measure $dH$ is an invariant
probability measure. Since $x$ is fixed, just compose $h$ with the
map that takes $\omega$ back to $\frac{x}{|x|}$ and conclude that
$h(\omega)=const.$

Going back we see that the expression in (\ref{stationary3}) equals
\begin{equation}
 \label{E:fourier3-1}
\int \widehat{\sigma}(|x| \xi) \, {|\widehat{\mu}(\xi)|}^2 \, d\xi
\end{equation}

\begin{equation}
 \label{E:fourier3-2}
= \lim_{\epsilon \to 0} \epsilon^{-1} \mu \times \mu \{(u,v): |x|
\leq |u-v| \leq |x|(1+\epsilon)\},
\end{equation}
so the problem reduces to showing that
\begin{equation}
 \label{E:mu-mu}
\mu \times \mu \{(u,v): |x| \leq |u-v| \leq |x|(1+\epsilon) \}
\lesssim \epsilon.
\end{equation}
This completes the proof of the three dimensional result, up to
(\ref{E:mu-mu}), which takes care of the first part of Theorem
\ref{main}. To prove the second part, observe again that we may
assume that $|x| \gtrsim 1$. By the method of stationary phase (see
e.g., \cite{St93}), we get
$$
|\widehat{\sigma}(\xi)| \lesssim {|\xi|}^{-1},
$$
and so the expression (\ref{E:fourier3-1}) is
$$
\lesssim \int {|\xi|}^{-1} {|\widehat{\mu}(\xi)|}^2 \, d\xi
=c \int \int {|x-y|}^{-2} d\mu(x) d\mu(y),
$$
which certainly converges if $\mu$ is the Frostman measure on a set
of Hausdorff dimension greater than two. This approach just fails to
work for two dimensional sets and this is where the $W^1_1({\Bbb
R}^2)$ assumption will play a key role. We now turn to the final
section of our paper where this is done. \vskip.125in

\section{Geometric estimates}

In this section we establish (\ref{E:mu-mu}) for measures supported
on graphs of $W^1_1({\Bbb R}^2)$ functions. Assume for a moment that
$\mu$ is the Lebesgue measure on a graph of a $C^1$ function $G$. We
may do that as long as our estimates do not quantitatively depend on
this smoothness assumption.

\begin{figure}
\begin{center}
\includegraphics[width=10cm]{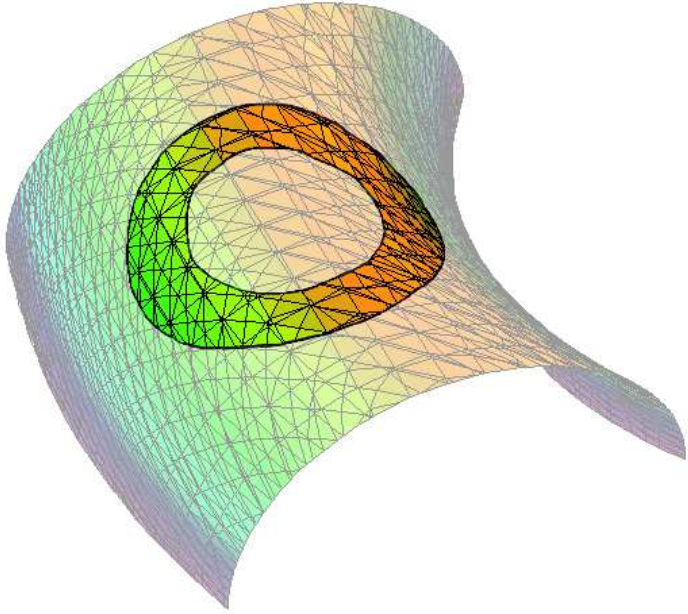}
\caption{ The set $S_x$ in the equation \eqref{E:Sx}. }
 \label{F:3D}
\end{center}
\end{figure}

We have
\begin{equation}
 \label{E:Sx}
S_x=\{y: 1 \leq |x-y| \leq 1+\epsilon\} \end{equation} is contained
in a curved annulus (see Figure \ref{F:3D}) whose dimensions are
$$
\epsilon \times \epsilon \times N,
$$
where
$$
N = \sup_{y,y' \in S_x} |G(y_1,y_2)-G(y'_1,y'_2)| \leq \epsilon
\sup_{z \in S_x} |\nabla G(z_1,z_2)|.
$$

If $G$ were Lipschitz, the right hand side above would automatically
be bounded by $\epsilon \, C_{Lip}$, independent of $x$, and the
proof of (\ref{E:mu-mu}) would be complete. Since we are only
assuming that $G$ is in $W^1_1({\Bbb R}^2)$, we have more work to
do. We must show that
$$ \int_{\{(x_1,x_2): x \in supp(\mu)\}}  \sup_{z \in S_x} |\nabla
G(z_1,z_2)| \,dx_1 \,dx_2 \leq C, $$
where $C$ does not depend on smoothness of $G$.

Since the set
$$
\{(z_1,z_2): z \in S_x\}
$$
is contained in a curved $\epsilon \times \epsilon$ square, and we
may take $\epsilon$ arbitrarily small, it is enough to show that
$$
\int_{\{(x_1,x_2): x \in supp(\mu)\}} |\nabla G(x_1,x_2)| dx_1dx_2 \leq C,
$$
and this follows instantly from the $W^1_1({\Bbb R}^2)$ assumption
on $G$. This completes the proof.

\end{document}